\documentclass[twoside,12pt,a4paper]{amsart} 
\usepackage{amsmath,amsfonts,amsthm,eurosym}
\usepackage[colorlinks=true,linkcolor=blue,citecolor=blue,filecolor=blue,urlcolor=blue,pageanchor=true,plainpages=false,linktocpage]{hyperref}
\usepackage[english]{babel}

\numberwithin{equation}{section}

\renewcommand{\epsilon}{\varepsilon}
\renewcommand{\phi}{\varphi}
\renewcommand{\rho}{\varrho}
\renewcommand{\theta}{\vartheta}

\begin{document}
\baselineskip3.6ex


\title{}

\author{Alexander Pigazzini$^{(*)}$}

\date{}

%


%

%

\makeatletter
\def\@makefnmark{}
\makeatother
\newcommand{\myfootnote}[2]{\footnote{\textbf{#1}: #2}}
 \footnote {}
 \footnote {}
{\bfseries\centerline{Curvature constrained on base for $(2+m)$-Einstein warped product manifolds}}
\\
\\
\centerline{A. Pigazzini, C. $\ddot{O}$zel and S. Jafari}
\\
\\
\section{abstract}
For the studied cases in \cite{Pigazzini}, the author showed that having the {\textit {$f$-curvature-Base}} ($R_{f_B}$) is equal to requiring a flat metric on the base-manifold.
In \cite{pndp} the authors introduced a new kind of Einstein warped product manifold, composed by positive-dimensional manifold and negative-dimensional manifold, the so called \textit{PNDP-manifolds}
The aim of this paper is to extend the work done in \cite{Pigazzini} to $m$-dimensional fiber showing if the value of $m$ can influence the result, i.e., finding base-manifolds with non-flat metric for $dimF \neq 2$, and doing some considerations of the $(2, m)$-PNDP manifolds with $R_{f_B}$.
\\
As a result, we find out that the dimension of fiber-manifold  does not change the result of \cite{Pigazzini}.
\\
Finally we add a \textit{Special Remark} about the possible use of the \textit{$(n,-n)$-PNDPs}, a special kind of Einstein warped product manifold, in superconductor graphene theory.
\\
\\
\textit{keywords}: \textit{$f$-curvature-Base}, $R_{f_B}$, Einstein warped product manifold, PNDP-manifolds.
\\
\\
\textit{2010 Mathematics Subject Classification}: 53C25, 53C21
\\
\section{Introduction and Preliminaries }
In recent years the study of \textit{warped product manifolds} (WPM) is of great interest both for the mathematicians and physicists. Many works have been published that have studied and introduced new types of WPM, (to name a few reference see \cite{Feitosa}, \cite{Tokura}, \cite{De} and \cite{Blaga}).
\\ Aytimur and Zgr in \cite{Aytimur} proved some results concerning the \textit{Einstein statistical WPM}, and in \cite{polj} Pigazzini et al. introduced a new type of \textit{WPM} so called \textit{PNDP-manifolds}, where the fiber is a manifold with negative dimension.
\\
In \cite{Pigazzini} Pigazzini introduced a simple constraint on the base-manifold called {\textit {$f$-curvature-Base}} ($R_{f_B}$) and proposed to use it in order to simplify the equations, trying to constructing a nonRicci-flat metric with non-constant Ricci curvature, on the base-manifold obtaining, as a result for the cases examined, that this is equivalent to the request for a flat metric. 
\\
\\
This paper is in effect an extension of the works done in \cite{Pigazzini} and moreover we make also a consideration about \cite{pndp}.
\\
\\
In the second part of the paper we reconsider the \textit{$(n,-n)$-PNDP}-manifolds and we suggest a possible use in the superconductor graphene theory.
\\
\\
{\bfseries Definition 2.1:}
A metric which satisfies the condition $Ric= \lambda g$ for some constant $\lambda$, is said to be an \textit{Einstein metric}. A manifold which admits an \textit{Einstein metric} is called an \textit{Einstein manifold}. (See \cite{Sambusetti}).
\\
\\
{\bfseries Definition 2.2:}
A warped product manifold is Einstein (see \cite{Pigazzini}, also \cite{Leandro}, \cite{Besse}) if and only if
\\
\numberwithin{equation}{section}
{(1)}
$\bar{Ric}=\lambda \bar{g} \Longleftrightarrow\begin{cases} 
Ric- \frac{m}{f}\nabla^2 f= \lambda g  \\  \ddot{Ric}=\mu \ddot{g} \\ f \Delta f+(m-1) |\nabla f|^2 + \lambda f^2 =\mu
\end{cases}$
\\
\\
where $\lambda$ and $\mu$ are constants, $m$ is the dimension of $F$, $\nabla ^2f$, $ \Delta f$ and $\nabla f$ are, 
\\
respectively, the Hessian, the Laplacian and the gradient of $f$ for $g$, with $f:(B) \rightarrow (0, \infty)$ a smooth positive function.
\\
Contracting first equation of (1) we get: 
\\
\numberwithin{equation}{section}
{(2)}
$R_Bf^2-mf \Delta f=n f^2 \lambda$ 
\\
where $n$ and $R_B$ is the dimension and the scalar curvature of $B$ respectively. By third equation, considering $m \neq 0$ and $m \neq 1$, we have:
\\
\numberwithin{equation}{section}
{(3)}
$mf\Delta f+m(m-1)|\nabla f|^2+m\lambda f^2=m\mu$
\\
Now from (2) and (3) we obtain:
\\
\numberwithin{equation}{section}
{(4)}
$|\nabla f|^2+[\frac{\lambda (m-n)+R_B}{m(m-1)}]f^2=\frac{\mu}{(m-1)}$
\\
\\
{\bfseries Definition 2.3:}
Let $(M, \bar{g})=(B,g)\times_f(F,\ddot{g})$ be an \textit{Einstein warped-product manifold} with $\bar{g}=g+f^2 \ddot{g}$. We define the scalar curvature of the Base-manifold $(B,g)$ as {\textit {$f$-curvature-Base}} ($R_{f_B}$), if it is a multiple of the warping function $f$ (i.e. $R_{f_B}=cf$ for $c$ an arbitrary constant belonging to $\mathbb{R}$). (See \cite{Pigazzini} as reference).
\\
\\
{\bfseries Remarks 2.1:} Since a warped product manifold (WPM) implies a non-constant warping function $f$ (otherwise it would be a simply Riemannian product-manifold, for more datails see \cite{Gebarowski}, \cite{O'Neill}), the results analyzed from now on will be considered from this point of view.
\\
\\
\\
{\bfseries {Case 1: Ricci-flat EWP ($\lambda=0$)} with nonRicci-flat fiber-manifold ($\mu \neq 0$)}
\\
\\
{\bfseries Theorem 2.1} \textit{Let $(M^{2+m}, \bar{g})=(B^2,g)\times_f(F^m,\ddot{g})$, be an \textit{Einstein warped-product \\ manifold} Ricci-flat (i.e., $\bar{Ric}=\lambda \bar{g}$ with $\lambda=0$), where $(B^2,g)$ is a smooth surface with non-zero $R_{f_B}$, and $(F^m,\ddot{g})$ is a smooth \textit{Einstein-surface} (i.e. $\ddot{Ric}=\mu \ddot{g}$).
\\
Then $(M^{2+m}, \bar{g})$, cannot exist}.
\\
\\
{\textit {Proof.}} In our case, we have $n=2$, $\lambda=0$ and $R_B=R_{f_B}$ (see \cite{Pigazzini}), then (2) and (3) become:
\\
\\
(5) $\Delta f - hf^2=0$\; \; \; (with $h=c/m$.)
\\
(6) $f\Delta f+(m-1)|\nabla f|^2+=\mu$
\\
Then (4) becomes:
\\
\numberwithin{equation}{section}
{(7)}
$(m-1)|\nabla f|^2+hf^3=\mu$
\\
\\
Now, by the initial hypothesis (non-zero $R_{f_B}$), we assume $h \neq 0$ with $f$ nonconstant and  setting $p=(m-1)$ and $u=-hf$, for an open set where $u$ nonzero. Thus:
\\
(8) $\Delta u + u^2=0$
\\
(9) $u\Delta u+p|\nabla u|^2- u^3 - h^2\mu =0$.
\\
(10) $p|\nabla u|^2- u^3 - h^2\mu =0$.
\\
\\
For the sake of semplicity we replace the constant $h^2 \mu$ with constant $A$.
\\
Let $g$ be the metric on $B$ and assume that $u$ is a nonzero (and hence necessarily positive) solution, to the above system on a simply-connected open subset $B' \subset B$.
\\
The equation (10) implies that $\omega_1=(u^3+A)^{-\frac{1}{2}} p^{\frac{1}{2}}du$, and this implies that we have to assume $(u^3+A)$ to be nonzero, $\omega_1$  is a $1$-form with $g$-norm $1$ on $B'$ and hence $g$ can be written in the form 
$g = \omega_1^2 + \omega_2^2$ for some $\omega_2$ which is also a unit $1$-form. 
\\
Fix an orientation by requiring that $\omega_1 \wedge \omega_2$ is the $g$-area form on $B'$, then 
\\
$\star du=(u^3+A)^{\frac{1}{2}} p^{-\frac{1}{2}} \omega_2$, and since $d(\star du)=\Delta u \; \omega_1 \wedge \omega_2$, it follows that:
\\
\\
$p^{-\frac{1}{2}}\frac{3}{2}(u^3+A)^{-\frac{1}{2}} u^2 du \wedge \omega_2 +p^{-\frac{1}{2}} (u^3 + A)^{\frac{1}{2}} d\omega_2=$
\\
$=d[(u^3+A)^{\frac{1}{2}} p^{-\frac{1}{2}} \omega_2]=-u^2 \omega_1 \wedge \omega_2=-u^2(u^3+A)^{-\frac{1}{2}} p^{\frac{1}{2}}du \wedge \omega_2=$
\\
$=\frac{3}{2}(u^3+A)^{-\frac{1}{2}}u^2du \wedge \omega_2+p u^2(u^3+A)^{-\frac{1}{2}}du \wedge \omega_2=-(u^3+A)^\frac{1}{2}d\omega_2$.
\\
\\
Then $(-\frac{3}{2}-p)(u^3+A)^{-1}u^2 du \wedge \omega_2=d\omega_2$ and we have $d[(u^3+A)^\frac{3+2p}{6} \omega_2]=0$,
\\ 
i.e. $\omega_2=(u^3+A)^{\frac{-3-2p}{6}} dv$, so the metric $g=p(u^3+A)^{-1}du^2+ (u^3+A)^{-1-\frac{2}{3}p}dv^2$ has a singularity in $u^3=-A$.
\\
\\
{\bfseries Remarks 2.2}: Here, the analysis of the singularity is substantially the same as in \cite{Pigazzini}. {\bfseries }
\\
\\
Consider an open set for which $(u^3+A)$ is nonzero. The Gaussian curvature is given by: $K=-\frac{9}{2}(u^3+A)^{-1}u^4p-\frac{9}{2}(u^3+A)^{-1}u^4p^2-(u^3+A)^{-1}u^4p^3+2up^3+3up^2$.
\\
In this case, it is easy to verify that for the initial hypothesis, where we have set $R_B=R_{f_B}$ (i.e. $K=-u\frac{m}{2}$), we observe that $K$ is incompatible with our analysis. In fact we have:
\\
\\
(11) $-u\frac{m}{2}=-\frac{9}{2}(u^3+A)^{-1}u^4p-\frac{9}{2}(u^3+A)^{-1}u^4p^2-(u^3+A)^{-1}u^4p^3+2up^3+3up^2$
\\
or
\\
$m=(u^3+A)^{-1}u^3(9p+9p^2+2p^3)-6p^2-4p^3$
\\
Now remembering that $m=p+1$ we have:
\\
$p+1+6p^2+4p^3=(u^3+A)^{-1}u^3(9p+9p^2+2p^3)$
\\
or
\\
$(p+1+6p^2+4p^3)(u^3+A)=u^3(9p+9p^2+2p^3)$
\\
or
\\
$(p+1+6p^2+4p^3)u^3+Ap+A+A6p^2+A4p^3=(9p+9p^2+2p^3)u^3$
\\
or
\\
$Ap+A+A6p^2+A4p^3=(-2p^3+3p^2+8p-1)u^3$
\\
Since $u$ must not be a constant, this implies:
\\
\\
(a) $4p^3+6p^2+p+1=0$ and (b) $-2p^3+3p^2+8p-1=0$. But the polynomials (a) and (b) have different solutions, so (11) is satisfied only for constant $u$ (i.e., $f=constant$), which is not admitted in our initial assumptions, therefore $(M^{2+m}, \bar{g})$, cannot exist.
\\
{\bfseries {Case 2: nonRicci-flat EWP ($\lambda \neq 0$)} with Ricci-flat fiber-manifold ($\mu = 0$)}
\\
\\
{\bfseries Theorem 2.2} \textit{Let $(M^{2+m}, \bar{g})=(B^2,g)\times_f(F^m,\ddot{g})$ be an \textit{Einstein warped-product \\ manifold}, where $(B^2,g)$ is a smooth surface with non-zero $R_{f_B}$, and $(F^m,\ddot{g})$ is a smooth Ricci-flat surface (i.e. $\ddot{Ric}=\mu \ddot{g}$, with $\mu=0$).
\\
Then $(M^{2+m}, \bar{g})$, cannot exist}.
\\
\\
{\textit {Proof.}} The analysis is essentially the same as seen so far, so we assume $h \neq 0$ and set $u=-hf$, where $f$ is not constant. The equations (2) and (4) become:
\\
\numberwithin{equation}{section}
{(12)}
$hf^2 - \Delta f - l \lambda f=0$
\\
\numberwithin{equation}{section}
{(13)} 
$|\nabla f|^2 + \frac{\lambda}{p}f^2 - \frac{l \lambda}{p}f^2+\frac{h}{p}f^3=0$
\\
where $h=\frac{c}{m}$, $l=\frac{2}{m}$ and $p=(m-1)$.
\\
Setting $u=-hf$ we obtain:
\\
\numberwithin{equation}{section}
{(14)}
$u^2+ \Delta u + Qu=0$
\\
\numberwithin{equation}{section}
{(15)} 
$|\nabla u|^2 -Su^2+Tu^2-Du^3=0$
\\
with $Q=\lambda l$, $S=\frac{\lambda l}{p}$, $T=\frac{\lambda}{p}$ and $D=\frac{1}{p}$, where it is easy to see that $D \neq 0$, $T \neq 0$, $S \neq 0$ and $S-T \neq 0$ for $m \neq 2$.
\\
\\
By the same token as in case (1b), we obtain from (15) that $du=(u^3D+Su^2-Tu^2)^{1/2}\omega_1$. This implies that we have to assume $(u^3D+Su^2-Tu^2)$ to be nonzero. \\ Then $\omega_1=(u^3D+Su^2-Tu^2)^{-1/2}du$, so $\star du=(u^3D+Su^2-Tu^2)^{1/2}\omega_2$.
\\
Since $d(\star du)=\Delta u \; \omega_1 \wedge \omega_2$, we obtain: 
\\ 
\\
\numberwithin{equation}{section}
{(16)} $d\omega_2=(\frac{-3}{2}u^2D-Su+Tu-u^2-Qu)(u^3D+Su^2-Tu^2)^{-1}du \wedge \omega_2$.
\\
\\
Thus we can write $\omega_2=u^{-A}(Du+S-T)^{-B}dv$ for some constant $A$ and $B$ and for some function $v$.
\\
Since for $u= \frac{T-S}{D}$ we have a singularity and we have assumed $(u^3D+Su^2-Tu^2)$ to be nonzero, then we must consider $u \neq \frac{T-S}{D}$ .
\\
\\
{\bfseries Remarks 2.3}:Even here, the analysis of the singularity is substantially the same as in the previous case (i.e. 1b, it is sufficient to consider $\frac{T-S}{D}=A$) so both are equivalent to the studied case in \cite{Pigazzini}.{\bfseries }
\\
\\
Continuing with the calculations we have: 
\\
$E=(u^3D+Su^2-Tu^2)^{-1}$ and
\\
$G=u^{-2A}(Du+S-T)^{-2B}$
\\
So by Brioschi's formula, we have that the Gaussian curvature is:
\\
$K=(2A^2+A)(Du+S-T)^3u^{(4+6A+4B)}+AD(Du+S-T)^2u^{(6A+4B+5)} \\ +(2ABD+4BD)u^{(6A+4B+7)}(Du+S-T)^{(2B+3)}\\ +(2B^2D^2+BD^2)u^{(6A+4B+8)}(Du+S-T)^{(2B+2)}$.
\\
Also in this case for the initial hypothesis, $2K=R_{f_B}=cf$, we must have $K=-u\frac{m}{2}$, which means: 
\\
\numberwithin{equation}{section}
{(17)}
$-\frac{m}{2}=(2A^2+A)(Du+S-T)^3u^{(3+6A+4B)}+AD(Du+S-T)^2u^{(6A+4B+4)} \\ +(2ABD+4BD)u^{(6A+4B+6)}(Du+S-T)^{(2B+3)}\\ +(2B^2D^2+BD^2)u^{(6A+4B+7)}(Du+S-T)^{(2B+2)}$.
\\
Now putting in relation the equation (16) with $\omega_2$, we obtain:
\\
\\
$\frac{-\frac{3}{2}u^2D-Su+Tu-u^2-Qu}{u^3D+Su^2-Tu^2}=\frac{-ADu-AS+AT-BDu}{Du^2+Su-Tu}$
\\
\\
and solving the partial fractions we have:
\\
\numberwithin{equation}{section}
{(18)} $A=\frac{Q}{Z}+1=\frac{m}{2-m}$, then $m \neq 2$ (and $m \neq 0$, $m \neq 1$ from \textit{Definition 2}).
\\
\numberwithin{equation}{section}
{(19)} $B=\frac{3D+2}{2D}-A=\frac{m-2m^2+2}{4-2m}$, then $m \neq 2$  (and $m \neq 0$, $m \neq 1$ from \textit{Definition 2}).
\\
where $Z=S-T$.
\\
\\
If (17) has a solution, certainly the coefficients of $u$ with highest degree must vanish. Hence we can consider the right side of (17) composed by:
\\
\\
$P_1(u)=(2A^2+A)(Du+S-T)^3u^{(3+6A+4B)}$ with highest degree: $6A+4B+6$,
\\
$P_2(u)=AD(Du+S-T)^2u^{(6A+4B+4)}$ with highest degree: $6A+4B+6$,
\\
$P_3(u)=(2ABD+4BD)u^{(6A+4B+6)}(Du+S-T)^{(2B+3)}$ with highest degree: $6A+6B+9$,
\\
$P_4(u)=(2B^2D^2+BD^2)u^{(6A+4B+7)}(Du+S-T)^{(2B+2)}$ with highest degree: $6A+6B+9$.
\\
\\
It is worth noticing that the highest degree of $P_1(u)$ is equal to that of $P_2(u)$ and the highest degree of $P_3(u)$ is equal to that of $P_4(u)$. But since the constants $A$ and $B$ can be non-integer and negative, we cannot know in advance which of the two degrees is the highest. We have 3 cases:
\\
\\
I) $6A+6B+9$ is the highest degree,
\\
II) $6A+4B+6$ is the highest degree,
\\
III) $6A+6B+9=6A+4B+6$.
\\
\\
\\
CASE (I):
\\
From coefficients of $P_3(u)$ and $P_4(u)$ if the (17) is satisfied, we should get: \\ $2A+2B+5=0$ and considering the (18) and (19) we have: 
\\
\numberwithin{equation}{section}
{(20)} $-4m^2-4m+24=0$, i.e., $m=2$ that is not possible for (18) and (19), and $m=-3$.
\\
Now if the CASE (I) (the highest degree) vanish for $m=-3$, we must consider the other degrees and also they must vanish for $m=-3$, so we proceed to consider the degree of CASE(II).
\\
\\
CASE (II):
\\
If the (17) is satisfied by considering coefficients of $P_1(u)$ and $P_2(u)$, we get: \\ $(2A^2+A)D^3+AD^3=0$ and since $D$ is nonzero we can divide for $D^3$, then:
\\
\numberwithin{equation}{section}
{(21)} $A=-1$
\\
Considering the (18): $-2=0$ is not possible regardless of the value of $m$.
\\
CASE (III):
\\
The equality in the CASE (III) implies $B=-\frac{3}{2}$ and for the (19) this means:
\\
\\
\numberwithin{equation}{section}
{(22)} $-2m^2-2m+8=0$, that has no solution for the integer values of $m$.
\\
\\
We showed that (17) could be satisfied only for some constant value of function $u$ (i.e. $f$ constant), which is not admitted in our initial assumptions. Then, also in this case, $(M^{2+m}, \bar{g})$, cannot exist.
\\
\\
\\
{\bfseries {Case 3: $(2, m)$-PNDP-manifolds with $R_{f_B}$}}
\\
\\
{\bfseries Remarks 2.4} The $(2, m)$-PNDP manifolds with $R_{f_B}$, does not exist.
\\
\\
 From PNDP-manifold definition (see \cite{pndp}), for $(2, m)$ case, we know that $dim \widetilde B=dimB'=1$, so it is well known that $1$-dimensional manifolds are Ricci-flat.  
\\
From \cite{Besse} we know that for $1$-dimensional base with Ricci-flat fiber (i.e., $\mu=0$), exists an Einstein warped product manifold with $\lambda=m$ and $f=e^t$. 
\\
Now if we consider $R_B=R_{f_B}$, we should have $R_B=ce^t$, but also being well known that for a product manifold the Ricci curvature of the product equals the sum of the Ricci curvatures of each manifolds of the product (see \cite{Atceken}), we obtain that the $(2, m)$-PNDP manifolds is Ricci-flat, so the scalar curvature of the base-manifold can not be $ce^t$. 
\\
\\
As known, the \textit{PNDP-manifolds} are born from the study of the Einstein-warped product manifolds, for this reason the following session wants to be dedicated to a possible and important application of the latters.
\\
\section{Special Remarks about $(n, -n)$-PNDP manifolds in Superconductors Graphene Model}
First of all we recall and highlight that the purpose of the \textit{PNDP-manifolds} is precisely to present the point-like manifolds from a mathematical point of view, and introduce a type of manifold with a new kind of hidden dimensions.
\\
In \cite{Capozziello}, Capozziello et al. introduced the concept of the "point-like manifold" building superconductors with graphene, in particular they argue that superconductor graphene can be produced by molecules organized in point-like structures where sheets are constituted by $(N+1)$-dimensional manifold. Particles like electrons, photons and ‘‘effective gravitons’’ are string modes moving on this manifold. In fact, according to string theory, bosonic and fermionic fields like electrons, photons and gravitons are
particular ‘‘states’’ or ‘‘modes’’ of strings. In their important work, they show that at the beginning, there are point-like polygonal manifolds (with zero spatial dimension) in space which strings attaching them, where all interactions between strings on one manifold are the same and are concentrated on one point which manifold is located on it. They also attaching to show that by joining these manifolds, 1-dimensional polygonal manifolds are emerged on which gauge fields and gravitons live and so, these manifolds glued to each other build higher dimensional polygonal manifolds with various orders of gauge fields and curvatures. 
\\
In this context, we think that the \textit{$(n, -n)$-PNDP manifolds} could play an important role.
In fact \textit{$(n,-n)$-PNDP} appears as a point (point-like), because in general, from our interpretation see \cite{pndp}, it is a point (positive and negative dimensions hidel each other out and and the total dimension equals zero), but in special it is composed by two manifolds, $B$ and $F$ with nonzero dimensions, so for the first time we have an object that looks like a point (point-like), but with a geometric structure on which we can make calculations.
\\
\\
Back to the graphene superconductors model, our \textit{$(n,-n)$-PNDP manifold} consists of two manifolds with nonzero dimensions (one with $n$-dimension and one with $-n$-dimension, where these two manifolds can be thought as a result of intersection between other manifolds). Then we can consider these two manifolds as contained in a "$p$-dimensional BULK", but their warped product (which generates the \textit{$(n,-n)$-PNDP)} will create the point-like polygonal manifold, a point-like space-time as supposed in \cite{Capozziello}.
\\
The \textit{$(n,-n)$-PNDPs} can be considered as possible mathematical interpretation of  point-like manifolds, because they render, for the first time, this abstract concept as a coherent mathematical object.
\\
\\
\\
{\bfseries \centerline{Conclusion}}
\\
\\
We have observed that, the dimension of fiber-manifold does not infuence the result with respect to what obtained in \cite{Pigazzini}. Not even the construction of a PNDP-manifold is made possible for 2-dimensional base-manifold case with $R_{f_B}$.
In conclusion we try to observe a possible important application for the $(n;-n)$-PNDP manifolds in the context of the superconductor graphene theory.
\\
\\
\\

\par \bigskip
\
\\\textit{Alexander Pigazzini}: IT-Impresa srl, 20900 Monza, Italy, E-Mail: pigazzinialexander18@gmail.com
\\\textit{Cenap Ozel}: King Abdulaziz University, Department of Mathematics, 21589 Jeddah KSA, E-Mail: cenap.ozel@gmail.com
\\\textit{Saeid Jafari}: College of Vestsjaelland South, Herrestraede 11, 4200 Slagelse, Denmark, E-mail: jafaripersia@gmail.com


\begin{thebibliography}{1}
	\bibitem{Aytimur} H. Aytimur and C. Özgür, \textit{Einstein Statistical Warped Product Manifolds}, \text{Filomat} 32:11 (2018), pp. 3891–3897.
	\bibitem{Besse} A.L. Besse - \textit{Einstein Manifolds} - \text{Springer-Verlag} (1987)
	\bibitem{Blaga} A. M. Blaga and C.E. Hretcanu, \textit{Remarks on metallic warped product manifolds}, \text{Facta Universitatis series: Mathematics and Informatics} Vol. 33, 2 (2018), pp. 269-277.
	\bibitem{De} U. C. De, S. Shenawy and B. Ünal, - \textit{Sequential Warped Products: Curvature and Killing Vector Fields}, \text{arXiv:1506.06056 [math.DG]}.
          \bibitem{Capozziello} S. Capozziello, R. Pincak, E. N. Saridakis, \textit{Constructing superconductors by grapheneChern–Simons wormholes}, \text{Annals of Physics (Elsevier)} 390 (2018), pp. 303-333.
	\bibitem{Feitosa} F. Feitosa, A. A. Freitas Filho, J. N. Vieira Gomesand R. Pina, \textit{On the construction of gradient almost Ricci soliton warped product}, \text{Nonlinear Analysis} 161 (2017), pp. 30-43.
	\bibitem{Gebarowski} A. Gebarowski, \textit{On Einstein warped products}, \textit{Tensor (N.S.)}, 52 (1993), pp. 204-207.
	\bibitem{Leandro} B. Leandro, M. Lemes de Sousa and R. Pina - \textit{On the structure of Einstein warped product semi-Riemannian manifolds} -  \text{Journal of Integrable Systems} 3 (2018), pp. 1–11. 
	\bibitem{O'Neill} B. O'Neill, \textit{Semi-Riemannian Geometry with Applications to Relativity}, \text{NewYork, Academic Press} (1983).
	\bibitem{Pigazzini} A. Pigazzini, \textit{On the $(2+2)$-Einstein Warped Product Manifolds with $f$-curvature-Base.} \text{Italian journal of pure and applied mathematics}, 44 (2020) (to appear).
           \bibitem{pndp} A. Pigazzini, C. $\ddot{O}$zel, P. Linker, S. Jafari, \textit{On PNDP-manifolds}, \text{\href{https://www.researchgate.net/publication/341326043_On_PNDP-manifolds}{doi:10.13140/RG.2.2.23516.82568/14}}. 
	\bibitem{Sambusetti} A. Sambusetti, \textit{Einstein Manifolds and Obstructions to the Existence of Einstein Metrics} \text{Rendiconti di Matematica, Serie VII}, 18 (1998), pp. 131-149.
	\bibitem{Tokura} W. I. Tokura, L. Adriano, R. Pina and M. Barboza, - \textit{Gradient Estimates on Warped Product Gradient Almost Ricci Solitons}, \text{arXiv:1905.00068 [math.DG]}.
\bibitem{Atceken} M. Atceken and S. Keles - \textit{On the product Riemannian manifolds} - 
\textbf{} Differential Geometry - Dynamical Systems (2003) p. 1-8.
\end{thebibliography}
\end{document}